\documentclass[11pt]{amsart}
\usepackage{graphicx}
\usepackage{amsfonts, amsmath, amsfonts, amssymb}
\newtheorem{thm}{Theorem}[section]
\newtheorem{cor}[thm]{Corollary}

\newtheorem{prop}[thm]{Proposition}
\theoremstyle{definition}
\newtheorem{defn}[thm]{Definition}
\theoremstyle{remark}
\newtheorem{rem}[thm]{Remark}
\newtheorem{quest}[]{Question}

\newtheorem{conj}[]{Conjecture}

\begin{document}

\title[Fat Triangulations and  Differential Geometry]
{Fat Triangulations and  Differential Geometry}

\author[Emil Saucan]{Emil Saucan}

\address{Department of Mathematics, Technion, Haifa, Israel}

\email{semil@tx.technion.ac.il}

\thanks{Research partly supported 
by European Research Council under the European Community's Seventh Framework Programme (FP7/2007-2013) / ERC grant agreement n${\rm ^o}$ [203134].}
\date{\today}

\maketitle


\begin{abstract}
We study the differential geometric consequences of our previous result on the existence of fat triangulations, in conjunction with a result of Cheeger, M\"{u}ller and Schrader, regarding the convergence of Lipschitz-Killing curvatures of piecewise-flat approximations of smooth Riemannian manifolds. A further application to the existence of quasiconformal mappings between manifolds, as well as an extension of the triangulation result to the case of almost Riemannian manifolds, are also given. In addition, the notion of fatness of triangulations and its relation to metric curvature and to excess is explored. Moreover, applications of the main results, and in particular a purely metric approach to Regge calculus, are also investigated.
\end{abstract}


\section{Introduction and Main Result}

The purpose of this note is to explore some differential geometric and function theoretic implications of the following theorem:

\begin{thm}[\cite{s2}] \label{thm:eu}
Let $M^n$ be a connected, oriented $n$-dimensional manifold ($n \geq 2$), with boundary, having a finite
number of compact boundary components, and such that
\begin{enumerate}
\item $M^n$ is of class $\mathcal{C}^{r},\; 1 \leq r \leq \infty\,, \; n \geq 2$;
\item $M^n$ is a $PL$ manifold and $n \leq 4$;
\item $M^n$ is a topological manifold and $n \leq 3$.
\end{enumerate}
If the boundary components admit fat triangulations of fatness $\geq \varphi_0$, then there exist a global fat
triangulation of $M^n$.
\end{thm}

\begin{rem}
In fact, the conditions on the compactness and boundedness of the boundary components in the theorem above are too strong, as indicated by the results in \cite{s06b}, \cite{s08a}, and their role is to exclude certain ``pathological'' cases.
\end{rem}

Given the triangulation results of \cite{Mun}, \cite{cms} for manifolds without boundary, the following corollary\footnote{See also \cite{Ca34}.} follows immediately:

\begin{cor} [\cite{s2}] \label{cor:eu}
Let $M^n$ be a connected, oriented $n$-dimensional manifold ($n \geq 2$), without boundary or having a finite
number of compact boundary components, and such that
\begin{enumerate}
\item $M^n$ is of class $\mathcal{C}^{r},\; 1 \leq r \leq \infty\,, \; n \geq 2$;
\item $M^n$ is a $PL$ manifold and $n \leq 4$;
\item $M^n$ is a topological manifold and $n \leq 3$.
\end{enumerate}
Then $M^n$ admits a fat triangulation.
\end{cor}

Recall that  or {\it fat} triangulations (also called {\it thick} in some of the literature) are defined (in \cite{cms}) as
follows:

\begin{defn} \label{def:fatness}
Let $\tau \subset \mathbb{R}^n$ ; $0 \leq k \leq n$ be a $k$-dimensional simplex.
The {\it fatness}  $\varphi$ of $\tau$ is defined as being:
\begin{equation} \label{eq:fatness}
\varphi = \varphi(\tau) = \hspace{-0.3cm}\inf_{\hspace{0.4cm}\sigma
< \tau
\raisebox{-0.25cm}{\hspace{-0.9cm}\mbox{\scriptsize${\rm dim}\,\sigma=j$}}}\!\!\frac{Vol_j(\sigma)}{{\rm diam}^{j}\,\sigma}\;.
\end{equation}
The infimum is taken over all the faces of $\tau$, $\sigma < \tau$,
and ${\rm Vol}_{j}(\sigma)$ and ${\rm diam}\,\sigma$ stand for the Euclidian
$j$-volume and the diameter of $\sigma$ respectively. (If
${\rm dim}\,\sigma = 0$, then ${\rm Vol}_{j}(\sigma) = 1$, by convention.)
A simplex $\tau$ is $\varphi_0${\it-fat}, for some $\varphi_0 > 0$, if $\varphi(\tau) \geq \varphi_0$. A triangulation (of a submanifold of $\mathbb{R}^n$) $\mathcal{T} = \{ \sigma_i \}_{i\in \bf I}$ is
$\varphi_0${\it-fat} if all its simplices are $\varphi_0$-fat. A
triangulation $\mathcal{T} = \{ \sigma_i \}_{i\in \bf I }$ is {\it
fat} if there exists $\varphi_0 \geq 0$ such that all its
simplices are $\varphi_0${\it-fat}.
\end{defn}

The proof of Theorem \ref{thm:eu} (hence that of Corollary \ref{cor:eu}) given in \cite{s2} essentially rests upon, 
among others, on a method of Cheeger et al. \cite{cms} for ``fattening'' triangulations.\footnote{For a different proof, without making appeal to \cite{cms}, see \cite{s06}.}
However, the ultimate 
goal of this fattening procedure and, indeed, of the whole paper \cite{cms}, is to prove 
Theorem \ref {thm:cms} 
below, regarding the approximation of curvature measures in $PL$ approximations. It is, therefore, natural to wish to apply our result mentioned above, in light of Theorem \ref{thm:cms}, that is to exploit its differential geometric aspects (and applications).

\begin{thm}[\cite{cms}] \label{thm:cms}
Let $M^n$ be a compact Riemannian manifold, with our without boundary, and let $M^n_i$ be a a sequence
of fat piecewise-flat manifolds converging to $M^n$ in the Hausdorff metric.
Denote by $\mathfrak{R}$ and $\mathfrak{R}_i$ respectively, the
{\em Lipschitz-Killing curvatures} of $M^n$, $M^n_i$.
 Then $\mathfrak{R}_i \rightarrow \mathfrak{R}$ in the sense of measures.
\end{thm}

Recall that, for a Riemannian manifold $M^n$, the Lipschitz-Killing curvatures are defined as follows:
\begin{equation} \label{eq:Rj-smooth}
R^j(M^n) = \frac{1}{(n-j)!2^j\pi^{j/2}(j/2)!}\sum_{\pi \in
S_n}(-1)^{\epsilon(\pi)}\Omega_{\pi(1)\pi(2)}\wedge\cdots\wedge\Omega_{\pi(j-1)\pi(j)}\wedge\omega_{\pi(j+1)}\wedge
\end{equation}
\[
\wedge\cdots\wedge\omega_{\pi(n)}\,,
\]
where $\Omega_{\pi(j-1)\pi(j)}$ are the {\it curvature $2$-forms} and
$\omega_{kl}$ denote the {\it connection $1$-forms}, and they are interconnected(interrelated) by the structure equations:

\begin{equation}
 \left\{
         \begin{array}{ll}
         d\omega_k = -\sum_i\omega_{kl}\wedge\omega_l\,,\\
         d\omega_{kl} = -\sum_i\omega_{ki}\wedge\omega_{il} + \Omega_{kl}\,.
         \end{array}
 \right.
\end{equation}
where $\{\omega_k\}$ is the dual basis of $\{e_k\}$\,.

\begin{rem}
The low dimensional Lipschitz-Killing curvatures are, in fact, quite familiar:  $R^0 \equiv$ {\rm volume} and $R^2 \equiv$ {\rm scalar curvature}. Moreover, $R^n  \equiv$  {\it Gauss-Bonnet-Chern form}, (for $n = 2k$).
\end{rem}

\begin{rem}
The integral $\int_{M^n}R^j$ is also known as the {\it integrated mean curvature (of order $j$)}.
\end{rem}

In a similar manner (but technically slightly more complicated), one can define the associated boundary curvatures (or {\it mean curvatures}) $H^j$ which are curvature measures on $\partial M^n$:
Let $\{e_k\}_{1 \leq k \leq n}$ be an orthonormal frame for the tangent bundle $T_{M^n}$ of $M^n$, such that, along the boundary $\partial M^n$, $e_n$ coincides with the inward normal. Then, for any $2k+1 \leq j \leq n$, we define
\begin{equation}
H^j = \sum_k\Omega_{j,k}\,,
\end{equation}
where

\begin{equation}
\Omega_{j,k} = c_{j,k}\sum_{\pi \in S_{n-1}}(-1)^{\epsilon(\pi)}\Omega_{\pi(1)\pi(2)}\wedge\cdots\wedge\Omega_{\pi(2k-1)\pi(2k)}\wedge\omega_{\pi(2k+1),n}\wedge\cdots
\wedge\omega_{\pi(j-1),n},\omega_{\pi(j)}\wedge
\end{equation}
\[
\wedge\cdots\wedge\omega_{\pi(n-1)}\,,
\]
and where

\begin{equation}
c_{j,k} =  \left\{
         \begin{array}{ll}
         (-1)^k\left(2^j\pi^{\frac{j-1}{2}}k!\left(\frac{j-1}{2}-k\right)!(n-j)!\right)^{-1}  \,, & j = 2p+1\\
         (-1)^k\left(2^{k+\frac{j}{2}}\pi^{\frac{j}{2}}k!(j-2k-1)!(n-j)!\right)  \,, & j = 2p\,.
         \end{array}
 \right.
\end{equation}

These curvatures measures are normalized by requesting(imposing the condition) that:
\begin{equation}
\int_{T^{n-j}\times M^j}R^j + \int_{T^{n-j}\times \partial M^j}H^j = \chi(M^j){\rm Vol}T^{n-j}\,,
\end{equation}
for any flat $T^{n-j}$.

\begin{rem}
As is the case with the Lipschitz-Killing curvatures, the low dimensional boundary curvatures also have quite familiar interpretations: $H^1 \equiv$ {area boundary}, $H^2 \equiv$ {mean curvature for inward
normal}\footnote{as expected given the generic names for these $H^j$-s}, etc.
\end{rem}



\begin{rem}
One can fatly triangulate the smooth manifold $M^n$ and obtain the desired approximation results for curvatures using the intrinsic metric, not just $PL$ (Euclidean) approximations \cite{reg}, \cite{cms} -- see also Section \ref{sec:almostRiem} below, in particular Formulas (\ref{eq:approx1}), (\ref{eq:approx2}) and Remark \ref{rem:gen}.
\end{rem}

The quest (indeed, the need) for fat triangulations as a prerequisite for Theorem \ref{thm:cms}, should not be surprising, in view of the fact that
fat triangulations can be characterized as having dihedral angles bounded from below (see Proposition \ref{eq:fat-CMS} and the discussion following it) and in view of the following expression of the Lipschitz-Killing curvatures in terms of dihedral angles (see \cite{cms} for the proof):
\begin{equation} \label{eq:Rj}
R^j = \sum_{\sigma^{n-j}}\big\{1 - \chi\big(L(\sigma^j)\big) + \sum_{l}\measuredangle(\sigma^{n-j},\sigma^{n-j+i_1})\cdots\measuredangle(\sigma^{n-j+i_{1-1}},\sigma^{n-j+i_1})\cdot
\end{equation}
\[
\cdot\left[1
- \chi\big(L(\sigma^{n-j+i_l})\big)\right] \big\}{\rm Vol}(\sigma^{n-j})\,,
\]
where $L(\sigma^{j})$ denotes the ({\it spherical}) {\it link} of $\sigma^{j}$, and $\measuredangle(\sigma^i,\sigma^j)$ is the internal dihedral angle of $\sigma^i < \sigma^j$; $\measuredangle(\sigma^i,\sigma^j) =
{\rm Vol}\big(L(\sigma^i,\sigma^{j})$, where the volume is normalized such that ${\rm Vol}(\mathbb{S}^n) = 1$,
for any  $n$.\footnote{See \cite{cms} for further details.}
(Here $\chi$, ${\rm Vol}$ denote, as usual, the Euler characteristic and volume of
$\sigma^k$, respectively.)

\begin{rem}
Theorem \ref{thm:cms} is given in terms of the intrinsic geometry of $M^n$. A similar characterization of the curvature measures in terms of the extrinsic geometry (of embeddings in $R^n$) is given in \cite{fu}.
\end{rem}

\begin{rem} \label{rem:big-fat}
The condition that the triangulation necessarily becomes arbitrarily fine is, in fact, too strong if the manifold contains large flat regions. (A motivational example, widely noted and exploited in Computer Graphics, is that of a round cylinder in $\mathbb{R}^3$.)  We have noted in \cite{s08}, \cite{saz0} the need and possibility of a triangulation with variable density of vertices, adapting to curvature. (See also \cite{fu}, where the hypothesis that the mesh of the triangulation converges to zero is discarded).
\end{rem}

The immediate differential geometric consequence 
of Theorems \ref{thm:cms} and \ref{thm:eu}, as well as Corollary \ref{cor:eu}  -- where by ``immediate'' we mean here that it can be directly inferred by applying the methods of \cite{cms} --  is the following:

\begin{thm} \label{thm: main}
Let $N = N^{n-1}$, 
be a not necessarily connected manifold, such that $N =\partial M, M = M^n$, where $M^n$ is, topologically, as in the statement of Theorem \ref{thm:eu}.
\begin{enumerate}
\item If $M,N$ are $PL$ manifolds, then the Lipschitz-Killing curvature measures of $N$ can be extended to those of $M$. More precisely, there exist Lipschitz-Killing curvature measures $\mathfrak{R} = \{R^j\}$ on $\bar{M} = M \cup N$,
 such that $\mathfrak{R}|_N = \mathfrak{R}_N$ and $\mathfrak{R}|_M  = \mathfrak{R}_M$, except on a regular (arbitrarily small) neighbourhood of $N$,
 where $\mathfrak{R}_N$, $\mathfrak{R}_M$ denote the curvature measures of $N, M$ respectively.
\item If $M,N$ are smooth manifolds, then $\mathfrak{R}|_N = \mathfrak{R}_N$ and $\mathfrak{R}|_M  = \mathfrak{R}_M$ hold only in the sense of measures.
\end{enumerate}
\end{thm}


\begin{rem}
Recall that $R^j|_{\partial M^n} = H^j$ and, in the case of $PL$ manifolds, it represents the contribution of the $(n-j)$-dimensional simplices that belong to the boundary. (For an explicit formula, see any of the formulas (3.23), (3.38) or (3.39) of \cite{cms}.)
\end{rem}

\begin{rem} \label{rem:extension}
In a sense, the theorem above can be viewed, in view of the previous Remark, as the ``reverse'' of the result of \cite{cms}, Section 8, regarding the convergence of the boundary measures.

This extension of the curvature measures from the boundary of the manifold to its interior has applications in (Discrete) General Relativity (see Section \ref{subsec:wormholes}) but also in Graphics (see also Section \ref{subsec:graphics}). In this later case, since usually the so called ``volumetric'' data $D^3$ is embedded in $\mathbb{R}^3$, the only relevant curvatures are those of the observable surface, i.e. the boundary $D^2 = \partial D^3$. It is interesting (and important) to obtain a coherent geometry for the whole data set, while keeping the distortion minimal. Here, by ``minimal'', we mean both the technical term as well as the fact that, by Theorem \ref{thm: main} above, the region where any distortion exists can be made as small as desired.
\end{rem}

\begin{rem}
Incidentally, one notes that the result above implies that the original proof, due to  Allendoerfer and Weil \cite{AW} -- triangulation and embedding based -- of the generalized Gauss-Bonnet Theorem, can be slightly improved to show that, both for polyhedral ($PL$), as well as for Riemannian manifolds, the Gauss-Bonnet (or, rather, Allendoerfer-Weil) Formula holds not only globally, but also point-wise (but, as we already noted, only in the sense of measures).
\end{rem}


\subsection{A generalization: Almost Riemannian manifolds} \label{sec:almostRiem}
Theorems \ref{thm:eu} and \ref{thm: main} above admit generalizations (see also Remark \ref{rem:gen} below), of which we bring here a rather direct one.
We begin with 
the following definition (cf. \cite{Se0}):

\begin{defn} \label{def:aRs}
A metric space $(M,d)$ is called an {\it almost Riemannian} space iff
\begin{enumerate}
\item $M$ is a smooth manifold;
\item There exists a (smooth) Riemannian metric $g$ on $M$ and a constant $C_0 > 0$, such that, for any $x \in M$, there exists a neighbourhood $U(x)$, such that
\begin{equation} \label{eq:aRs}
C_0^{-1}d(y,z) \leq {\rm dist}_g(y,z)  \leq  C_0d(y,z)\,,
\end{equation}
for all $y,z \in U(x)$.
\end{enumerate}
\end{defn}

The basic example  of an almost metric space (beyond the trivial one $d \equiv g$) is given by any smooth (Riemannian) manifold embedded in some $\mathbb{R}^N$ and $d$ be the Euclidean distance in $\mathbb{R}^N$, $d = {\rm dist}_{Eucl}$, i.e. precisely the setting which we are concerned: the secant approximation of an embedded smooth manifold, with its Euclidean (ambient) metric is a almost Riemannian manifold (relative, so to say, to the approximated smooth manifold). According to \cite{Se0}, this holds, in fact, for any $C_0 > 0$, albeit at the price of locality (i.e. the size of the neighbourhood $U(x)$ depends very much on the point $x$) and cannot be apriorily be asserted (and, in fact, on the large scale, the to distance (geometries) might -- and usually do -- differ quite widely). However, we can do a bit better, by requiring that the simplices of the approximations are fat. Even if this requires some ``preprocessing'' of the given manifold this can be done (this being the starting point of this paper). Then, for any triangulation patch, we have the following estimates:
\begin{enumerate}
\item If $M$ has no boundary, then, by \cite{Pe}
\begin{equation} \label{eq:approx1}
\frac{3}{4}d_g(y,z) \leq d_{Eucl}(y,z) \leq  \frac{5}{3}d_g(y,z)\,;
\end{equation}
%
\item If $M$ has boundary, then, by \cite{s2}
\begin{equation} \label{eq:approx2}
\frac{3}{4}d_g(y,z) - f(\theta)\eta_\partial \leq d_{Eucl}(y,z) \leq  \frac{5}{3}d_g(y,z)  + f(\theta)\eta_{\partial M}\,;
\end{equation}
where $f(\theta)$ is a constant depending on the $\theta =
\min{\{\theta_{\partial M},\theta_{{\rm{int}}\, M}\}}$ -- the fatness
of the triangulation of $\partial M$ and ${\rm{int}\, M}$,
respectively, and $\eta_\partial$ denotes the {\it
mesh} of the triangulation
\end{enumerate}

Even though in the general case we can not produce estimates as precise as (\ref{eq:approx1}) and (\ref{eq:approx2}) above, (\ref{eq:aRs}) still holds.  
Since $M^n$ is a topological manifolds, we can triangulate it (using Munkres' classical results) and, furthermore, we can pass from $d(y,z)$ to ${\rm dist}_g$ and back in a controlled manner, using (\ref{eq:aRs}), thus allowing us to apply the  fattening techniques of \cite{cms}.\footnote{The $\varepsilon$-{\it moves} employed in \cite{cms} show that only the metric properties suffice to ensure fatness -- see Section \ref{sec:def-fat1} below for an extended discussion.}
It follows that Theorems \ref{thm:eu} 
generalizes to almost Riemannian manifolds, and
we formalize this observation as

\begin{thm} \label{thm:aRs}
Let $(M,d)$ be an almost Riemannian manifold, where $M$ satisfies the conditions in the statement of Theorem \ref{thm:eu}.
Then it admits a fat triangulation.
\end{thm}

\begin{rem}
Of course, one would also like 
to obtain a version of Theorem \ref{thm: main} adapted to almost Riemannian manifolds. However, this eludes us so far, since even a proper definition of the Lipschitz-Killing curvatures of an almost Riemannian manifold is not quite evident.
\end{rem}

\begin{rem} \label{rem:gen}
Besides the slight generalization (and its applications) presented above, there exist other, perhaps more less immediate ones -- see \cite{s08a}, \cite{Sa11}. In particular, we noted the existence of fat triangulations of Lipschitz manifolds. Since Alexandrov spaces are, by results of Perelman  \cite{Per} and Otsu and Shioya \cite{OS}, Lipschitz manifolds a.e., fat triangulations of such objects, and their differential geometric implications are certainly worth to explore, especially the role of Alexandrov spaces in Graphics, Imaging, etc., and also in Regge calculus -- see Section 3.1 below and also \cite{Sa11b}.

Note also that Semmes stresses that the condition in the definition of almost Riemannian manifolds, of $M$ being a topological manifold is far too strong. We can still use Semmes' original Definition \ref{def:aRs} for $PL$ (or polyhedral) manifolds (at least in dimension $n \leq 4$) by considering {\it smoothings} (see \cite{Mun}) and also considering, instead of the original $PL$ (polyhedral) metric $g$, its smoothing $\tilde{g}$ -- see \cite{Sa11b} for details and another application of this technique. (In fact, any topological manifold of dimension $n \leq 3$, admits a $PL$ structure -- see, e.g., \cite{Mun} -- therefore one can construct almost Riemannian manifolds of dimension 2 and 3, by starting with a topological manifold and considering, for example, the combinatorial metric on a compatible $PL$ structure, etc.)
\end{rem}

\begin{rem}
While, as we have noted above, not all the differential geometric consequences of the existence of fat triangulations of almost Riemannian manifolds are, at this point in time, accessible 
to us, this generalization is not gratuitous, either. Indeed, the fact that $PL$ manifolds -- satisfying certain additional technical conditions (see \cite{Se0}), that we do not bring here in order not to deviate to much from the main ideas -- are almost Riemannian manifolds, shows that they also admit differential forms satisfying additional properties -- see \cite{Se0}, Theorems 17.3 and 17.10. This differential forms are related to the construction of certain Lipschitz mappings that are connected, in their turn, to the existence of sufficiently large families of curves connecting two given points, 
with lengths not deviating too much 
from the distance between the two points \cite{Se0}. We shall investigate elsewhere the implications (and relevance) of the existence of such curves for $PL$ manifolds \cite{Sa11c}. However, we formalize the fact above as:

\begin{prop}
Let $(M,d)$ be a $PL$ manifold. Then there exists on $M$ bounded measurable (resp. locally integrable) differential forms that satisfy the conditions in \cite{Se0}, Theorem 17.3 (resp. \cite{Se0}, Theorem 1.10).
\end{prop}

Moreover, given the following facts: (a) If $\{P_m\}_{m \in \mathbb{N}}$ is a sequence of $n$-dimensional polyhedral manifolds, converging, in secant approximation to $M^n \subset \mathbb{R}^N$, then the convergence is also in the Lipschitz sense (see \cite{fu}, Section 3); (b) Theorems Theorem 17.3 and 1.10 of \cite{Se0} are proved using solely properties of Lipschitz functions (with values in $\mathbb{S}^n$), we formulate the following conjecture, whose full proof we shall bring elsewhere \cite{Sa11c}:

\begin{conj}
Let $M^n \subset \mathbb{R}^N$ be a $C^{1,1}$ compact manifold with boundary, and let $\{P_m\}_{m \in \mathbb{N}}$ be a sequence of fat $PL$ (polyhedral) manifolds with boundary {\it closely inscribed}\footnote{See \cite{fu}, p. 179 for the precise definitions.} in $M^n$, Hausdorff converging to $M^n$. Denote by $\theta_{0,m}, \theta_{T,m}, \Lambda_{1,m}$ and $\Lambda_{2,m}$ the bounded measurable, respectively locally integrable, differential forms given 
by \cite{Se0}, Theorems 17.3 and 1.10, respectively. Then $\theta_{0,m} \rightarrow \theta_0, \theta_{T,m}  \rightarrow \theta_T, \Lambda_{1,m} \rightarrow \Lambda_1, \Lambda_{2,m} \rightarrow \Lambda_2$,
where $\theta_0,\theta_T, \Lambda_1, \Lambda_2$ are bounded measurable, respectively locally integrable differential forms on $M^n$, satisfying the conditions in \cite{Se0}, Theorems 17.3 and 1.10, respectively. Here the convergence should be understood in the sense that it occurs only where it makes sense, that is on the space of forms on the common set of $\{P_m\}_{m \in \mathbb{N}}$ and $M^n$, i.e. on the vertices of the polyhedral manifolds $P_m$.
\end{conj}

\end{rem}


\section{On the definition of fat triangulations} \label{sec:def-fat}

\subsection{Comparison of the various definitions of fatness} \label{sec:def-fat1}

As already noted, the definition of fatness given above is that 
introduced in \cite{cms}. One reason for doing this is to preserve the ``unity of style'', so to say: since we heavily rely, at least in the first part, on the results and techniques of \cite{cms}, we find only fitting that we used, at least in the beginning, the same definition as that of Cheeger et al. However, they also prove in the same paper that the following result holds:

\begin{prop}[\cite{cms}] \label{eq:fat-CMS}
There exists a constant $c(k)$ that depends solely upon the
dimension $k$ of  $\tau$ such that
\begin{equation}   \label{eq:ang-cond}
\frac{1}{c(k)}\cdot \varphi(\tau) \leq \min_{\hspace{0.1cm}\sigma
< \tau}\measuredangle(\tau,\sigma) \leq c(k)\cdot \varphi(\tau)\,,
\end{equation}
and
\begin{equation}    \label{eq:area-cond}
\varphi(\tau) \leq \frac{Vol_j(\sigma)}{diam^{j}\,\sigma} \leq c(k)\cdot \varphi(\tau)\,,
\end{equation}
where $\varphi$ denotes the fatness of the simplex $\tau$, $\measuredangle(\tau,\sigma)$ denotes the  ({\em
internal}) {\em dihedral angle} of the face $\sigma < \tau$ and $Vol_{j}(\sigma)$; $diam\,\sigma$ stand for the
Euclidian $j$-volume and the diameter of $\sigma$, respectively. (If $dim\,\sigma = 0$, then $Vol_{j}(\sigma) = 1$,
by convention.)
\end{prop}

Note first that condition  (\ref{eq:ang-cond}) is just the expression of fatness as
a function of dihedral angles in all dimensions. This fact warrants a number of remarks:
\begin{itemize}
\item Using the dihedral angle approach in assuring the ``aspect ratio'' of a triangular mesh is commonly used in Computational Geometry, Computer Graphics and related fields, to ensure that the constituting tetrahedra are not to ``flat'' or too ``slim'', in particular that no ``slivers''' appear (see, e.g. \cite{AB}, \cite{D++}, \cite{Ed}).
\item The condition under scrutiny shows that thickness is hierarchical, in the sense that for a simplex to be thick, all its lower dimensional faces have to be thick. This is also transparent from condition \ref{eq:area-cond} and we shall discuss this fact again shortly, from a different point of view.
\item The dihedral angle definition appears to be quite promising  in the quite recent developments of Regge calculus (see, e.g. \cite{BD}, \cite{DS}, \cite{RW}). Indeed, it seems that, from the Theoretical Physics, it has quite a number of advantages. However, it falls short from Regge's original goal, as stated in \cite{reg}, for a purely {\it metric} gravity (both classical and quantum).
    Moreover, from the mathematicians point of view it also lacks the kind of ``symmetry'' one usually strives for, since it makes appeal to both angles and distances (e.g. in the ``fattening'' technique (using $\varepsilon$-{\it moves}) -- see \cite{cms}). More importantly, one strives for the most general possible setting, thus one looks for the possible extensions of the results of \cite{cms} (and, in consequence, of the present paper) to as ``general'' metric spaces. Therefore, a purely metric approach for all aspects of the problem, including fatness, is desirable and we shall discuss it in detail below aa well as in Section 3.1.2.
\end{itemize}

Similarly, condition (\ref{eq:area-cond}) expresses fatness as given by ``large
area/diameter'' (or ``volume/diameter'') ratio. Diameter is important since fatness is
independent of scale. Note, again, the hierarchic property 
of fatness as expressed by this condition.

Moreover, even the defining condition of fatness (\ref{eq:fatness}) lacks this purely metric aspect. Indeed, while volume and edge lengths of a Euclidean simplex (or, for that matter, in any space form) are closely related, in a general metric space no volume is apriorily postulated. 
Admittedly, one can attempt a fitting definition for {\it metric measure spaces} \cite{Sa11}, but this has only limited purely mathematical applications and, furthermore, it also departs from the purely metric approach.\footnote{even if such an approach is, in the light of the literature mentioned above, less useful in Physics than previously believed} Fortunately, this problem is easy to amend, using the so called {\it Cayley-Menger determinant}, that expresses the volume of the $n$-dimensional Euclidean simplex $\sigma_n(p_0,p_1,\ldots,p_n)$ of vertices $p_0,p_1,\ldots,p_n$ as a function of its edges $d_{ij}, 0 \leq i < j \leq n$:
\begin{equation} \label{eq:CayleyMenger-nD}
 \Gamma(p_0,p_1,\ldots,p_n) = \left| \begin{array}{ccccc}
                                            0 & 1 & 1 & \cdots & 1 \\
                                            1 & 0 & d_{01}^{2} & \cdots & d_{1n}^{2} \\
                                            1 & d_{10}^{2} & 0 & \cdots & d_{1n}^{2} \\
                                            \vdots & \vdots & \vdots & \ddots & \vdots \\
                                            1 & d_{n0}^{2} & d_{n1}^{2} & \cdots & 0
                                      \end{array}
                               \right| \,;
\end{equation}
namely
\begin{equation} \label{eq:VolEucl-n-simplex}
{\rm Vol}^2(\sigma_n(p_0,p_1,\ldots,p_n)) = \frac{(-1)^{n+1}}{2^n(n!)^2}\Gamma(p_0,p_1,\ldots,p_n)\,.
\end{equation}
(Similar expressions for the Hyperbolic and Spherical simplex also exist, see, e.g. \cite{Bl}, \cite{BM}.)

Evidently, that the Cayley-Menger determinant makes sense for any {\it metric $(n+1)$-tuple}, i.e. for any metric metric space with $n+1$ points $p_0,p_1,\ldots,p_n$  and mutual distances $d_{ij}, 0 \leq i \leq j \leq n$.

As implicitly mentioned above, fatnesss' Definition \ref{eq:fatness} is only one of the (many) existing ones. An alternative one, given solely in terms of distances, is given in \cite{Tu}.  Tukia's fatness\footnote{In fact he defines the reciprocal quantity which he calls {\it flatness} and denotes by $F(\sigma)$.} of a simplex $\sigma = \sigma(p_0,\ldots,p_n)$ can be written as
\begin{equation}
\varphi_T(\sigma) = \max_{\pi}{\frac{\delta(p_0,\ldots,p_n)}{{\rm diam}{\sigma}}}\;,
\end{equation}
where $\delta = \delta(p_0,\ldots,p_n) = {\rm dist}(p_i,F_i)$, where $F_i$ denotes the $(n-1)$-dimensional face (of $\sigma$) opposite to $p_i$, and where the maximum is taken over all the permutations $\pi$ of $0,\ldots,n$.
In other words, $\delta$ represents the maximal height $h_{Max}$ from the vertices of $\sigma$ to its $(n-1)$-dimensional faces (or {\it facets}).
Note that this definition of fatness can also be used, for instance, for Hyperbolic simplices (as, indeed, Tukia did).

Since the diameter of a simplex is nothing but the longest edge $l_{Max}$, that is ${\rm diam}{\sigma} = \max_{0 \leq i<j \leq n}d_{ij} = \max_{1 \leq m \leq n(n+1)/2}l_m$, Tukia's definition of fatness is, clearly, purely metric. It is true that the expression of $\varphi_T$ is not given solely in terms of the distances between the vertices (of $\sigma$), however it is easy to remedy this by computing $h_i$ -- the distance from $p_i$ to the face $F_i$, using the classical and well known formula:
\begin{equation} \label{eq:volume-height}
{\rm Vol}(\sigma) = \frac{1}{n}h_i{\rm Area}(F_i)
\end{equation}
and then expressing the volume of $F_i$ using the fitting Cayley-Menger determinant, where we denoted the $(n-1)$-volume by ``${\rm Area}$''.
(One has to proceed somewhat more carefully for the case, say, of Hyperbolic simplices, but we are not concerned here with this case.)

Peltonen's definition is, perhaps the easiest to express, even not if the simplest for actual computation, as she defines fatness (of a simplex) as:
\begin{equation}
\varphi_P = \frac{r}{R}\,,
\end{equation}
where $r$ denotes the radius of the inscribed sphere of $\sigma$ ({\it inradius})
and R denotes the radius of the circumscribed sphere of
$\sigma$ ({\it circumradius}).
The problem with this approach is that it is not given -- explicitly, that is -- in terms of the lengths of the edges of the simplex. Again, this does not represent a serious inconvenience, since for $R$ we have the following formula:
\begin{equation}
R = -\frac{1}{2}\frac{\Delta(p_0,\ldots,p_n)}{\Gamma(p_0,\ldots,p_n)}\,,
\end{equation}
where $\Delta(p_0,\ldots,p_n)$ is the minor of the element in the first line and first row of $\Gamma(p_0,\ldots,p_n)$, where
\begin{equation} \label{eq:CayleyMenger-nD}
 \Gamma(p_0,\ldots,p_n) = \left| \begin{array}{cccc}
                                            0 & d_{01}^{2} & \cdots & d_{1n}^{2} \\
                                            d_{10}^{2} & 0 & \cdots & d_{1n}^{2} \\
                                            \vdots & \vdots & \ddots & \vdots \\
                                            d_{n0}^{2} & d_{n1}^{2} & \cdots & 0
                                      \end{array}
                               \right| \,.
\end{equation}

Moreover, $r$ can be computed quite simply using the fact (akin to formula (\ref{eq:volume-height})) that
\begin{equation}
{\rm Vol}\sigma_n = \frac{1}{n}\sum_{1}^{n+1}{\rm Area}(F_i)\,,
\end{equation}
and expressing again ${\rm Area}(F_i)$ as a function of its edges.

On the other hand, the equivalence of Tukia's and Peltonen's definitions follows using Peltonen's computation of the $h_i$-s (and ${\rm diam}{\sigma}$), which is given -- {\it inter alia} -- in \cite{Pe}. However, the proof requires additional notations and definitions as well as some quite extensive technical details. Therefore, we refer the reader to the original paper of Peltonen where the connection between $\varphi_P$ and $h$ (hence $\varphi, \varphi_T$) is given (even if not quite explicitly).


Munkres' definition is, for the specific case of Euclidean simplices\footnote{But only in this case, as it represents a generalization of the Euclidean case for general simplices (as befitting Differential Topology goals).}, somewhat of a 
compromise between Cheeger's and Tukia's ones:
\begin{equation} \label{eq:fat-Minkres}
\varphi_M = \frac{{\rm dist}(b,\partial\sigma)}{{\rm diam}\,\sigma}\,,
\end{equation}
where $b$ denotes the {\it barycenter} of $\sigma$ and $\partial\sigma$ represents the standard notation for the {\it boundary} of $\sigma$ (i.e the union of the $(n-1)$-dimensional faces of $\sigma$). From the considerations above (and from \cite{Mun}, Section 9) it follows that for the Euclidean case, this definition of fatness is also equivalent to the previous ones.

\begin{rem} \label{rem:Fu-fat}
Fu \cite{fu} also introduces a definition of fatness that (up to the quite different notation) is identical to that of \cite{cms}, when restricted to individual simplices. However, for triangulations, his definition exceeds, in general, 
the one given in \cite{cms}, due to the fact that Fu discards, in his use of fatness, the requirement that the approximations become arbitrarily fine (see also Remark \ref{rem:big-fat} above).
\end{rem}

\begin{rem} 
Of course, there exists a (well know) duality between the spherical distances and dihedral angles, as expressed by the Gramm determinant (see, e.g. \cite{Ko} and the references therein), but unfortunately this duality does not hold precisely in the case that we are interested in, that is constant sectional curvature $K \equiv 0$, so we do not discuss it here. Let us note only the fact that this type of analysis involves a generalization of the Cayley-Menger determinant.
\end{rem}

To complete the circle -- so to say -- we should emphasize that, given the dihedral angles and the areas of the $(n-1)$-faces of the Euclidean (piecewise flat) triangulation, one can calculate the edge lengths. For example, one can devise explicit computations for the case important in Ricci calculus, i.e. 4-simplices, minus a number of so called ``dangerous configurations'', (and implicit ones for the general case) in \cite{DS}. (A related formula -- explicit for the 3-dimensional case -- this time involving the Cayley-Menger determinant, can be found in \cite{BM}, p. 344.)


\subsection{Fatness, metric 
curvature and excess}

Triangle thickness is related, quite directly, to another (geo-)metric invariant -- except curvature -- and, through it, again to curvature, both metric and, in the classical case of Riemannian manifolds, to Ricci curvature.\footnote{This connection is mainly through the Abresch-Gromoll theorem and its applications -- see, e.g. \cite{Zh} and, for a brief overview, \cite{Be}.} The quantity in question is the so called {\it excess}:

\begin{defn}
Given a triangle\footnote{not necessarily geodesic} $T = \triangle(pxq)$ in a metric space $(X,d)$, the {\it excess} of $T$ is defined as
\begin{equation}
e = e(T) = d(p,x) + d(x,q) - d(p,q).
\end{equation}
\end{defn}

The connection with thickness is obvious: a 2-simplex is fat iff it its excess is bounded away from 0. Note, however, that one can not assert that higher dimensional simplices are thick iff all their 2-dimensional faces are. Indeed, one can easily construct, for example, 3-dimensional tetrahedra with regular base and the other faces congruent to each other and having the common vertex at $\varepsilon$ distance from the barycenter of the base.

A local version of this notion -- introduced, it seems by Otsu \cite{Ot} -- also exist, namely the {\it local excess} (or, more precisely, the {\it local $d$-excess}):
\begin{equation} \label{eq:local-excess}
e(x) = \max_p\max_{x \in B(p,d)}\min_{q \in S(p,d)}\left(e(\triangle(pxq)\right)\,,
\end{equation}
where $d \leq {\rm rad}(X) = \min_p\max_qd(p,q)$, (and where $B(p,d), S(p,d)$ stand -- as they commonly do -- for the ball and respectively sphere of center $p$ and radius $d$).

Global variation of this quantity have also been defined:
\begin{equation}
e(X) = \min_{(p,q)}\max_x\left(e(\triangle(pxq)\right)\,,
\end{equation}
and, the so called\footnote{after \cite{Ot}} {\it global big excess}:
\begin{equation}
E(X) = \max_q\min_p\max_x\left(e(\triangle(pxq)\right)\,.
\end{equation}

\begin{rem}
The geometric ``content'' of the notion of local excess is that, for any $x \in B(p,d)$, there exists a (minimal) geodesic $\gamma$ from $p$ to $S(p,d)$ such that $\gamma$ is close to $x$.
\end{rem}

The relevant notion of curvature, more precisely a {\it metric} curvature, is the so called {\it Finsler-Haantjes curvature}\footnote{Introduced by Haantjes \cite{Ha}, based upon an idea of Finsler.}

\begin{defn}
 Let (M,d) be a (geodesic) metric space and let $c: I=[0,1] \stackrel{\sim}{\rightarrow} M$ be a homeomorphism, and let $p,q,r \in c(I),\; q,r \neq p$. Denote by
$\widehat{qr}$ the arc of $c(I)$ between $q$ and $r$. 
 Then $c$ has {\em Finsler-Haantjes curvature} $\kappa_{FH}(p)$ at the point $p$ iff:
\begin{equation}
\kappa_{FH}^2(p) = 24\lim_{q,r \rightarrow
p}\frac{l(\widehat{qr})-d(q,r)}{\big(l(\hat{qr})\big)^3}\,,
\end{equation}
where ``$l(\widehat{qr})$'' denotes the length -- in intrinsic
metric induced by $d$ -- of $\widehat{qr}$.

\end{defn}

\begin{rem}
A number of observations are mandatory: 
\\
(a) The Finsler-Haantjes curvature is obviously defined only for rectifiable curves, but this does not represent an impediment for our purposes, since we are concerned only with such curves.
\\
(b) As a curvature for curves $\kappa_{FH}$ is, fittingly, unsigned.
\\
(c) The constant 24 (or, rather, $\sqrt{24}$) in the definition of $\kappa_{FH}$ is not truly important for the definition of Finsler-Haantjes curvature in the most general possible context. Its role is to equate $\kappa_{FH}$ with the classical notion as defined for smooth curves in $\mathbb{R}^2$ ($\mathbb{R}^n$).
\\
(d) Alternatively, since for points/arcs where Finsler-Haantjes exists, \\ $\lim_{d(q,r) \rightarrow 0}\frac{l(\widehat{qr})}{d(q,r)} = 1$ (see \cite{Ha}),  $\kappa_{FH}$ can be defined (see, e.g. \cite{Ka}) by
\begin{equation}
\kappa_{FH}^2(p) = 24\lim_{q,r \rightarrow
p}\frac{l(\widehat{qr})-d(q,r)}{\big(d(q,r))\big)^3}\,\,;
\end{equation}
\end{rem}

The connection between thickness and Finsler-Haantjes curvature has its counterpart in a (quite expected) curvature-excess connection:

Namely, using a simplified notation and discarding (for convenience/simplicity) the normalizing constant, one has the following relation between the two notions:
\begin{equation}
\kappa_{FH}^2(T) = \frac{e}{d^3}\,,
\end{equation}
where by the curvature of a triangle $T = T(pxq)$ we mean the curvature of the path $\widehat{pxq}$.
Thus Finsler-Haantjes curvature can be viewed as a {\it scaled} version of excess. Keeping this in mind, one can define also a global version of this type of metric curvature, namely by defining, for instance:
\begin{equation}
\kappa_{FH}^2(X) = \frac{E(X)}{{\rm diam}^3(X)}\,,
\end{equation}
or
\begin{equation}
\kappa_{FH}^2(X) = \frac{e(X)}{{\rm diam}^3(X)}\,,
\end{equation}
as preferred.

Of course, one can proceed in the ``opposite direction'', so to speak, and express the proper (i.e. point-wise) Finsler-Haantjes curvature via the 
definition (\ref{eq:local-excess}) of local excess, as
\begin{equation}
\kappa_{FH}^2(x) = \lim_{d \rightarrow 0}e(x)\,.
\end{equation}

We do not explore here the full depth of the interconnections between these three concepts, and we postpone such analysis 
for further research (see \cite{Sa11c}), except the case of principal curvatures which we shall 
investigate in some detail in the 
next section. However, it  is important to underline the fact that, to obtain good approximations of curvatures of $PL$ manifolds, one needs, in fact, to start with 
triangulations having principal metric curvatures -- see discussion in Section \ref{sec:applic} below. (The other requirement is, of course, the metric approximation being arbitrarily fine, that is that the mesh of triangulation converges to zero.)

\begin{rem}
Amongst the different types of metric curvatures for curves, the {\it Menger curvature} (see, e.g. \cite{Bl}, \cite{BM}) is the best known and it has been applied quite frequently and successfully, not only (and, we could add, quite naturally) in curve reconstruction (see \cite{Gi}), but also, more deeply (and more interestingly), in computing estimates (obtained via the Cauchy integral) for the regularity of fractals and the fatness of sets in
the plane (see \cite{Pa}). While equivalent to the Finsler-Haantjes curvature in most metric spaces of interest (see \cite{Bl}, \cite{BM}), and certainly in our context, we have preferred here the later, due to transparency of its connection with the notion of excess, thence to fatness, as detailed above.
\end{rem}


\section{Applications} \label{sec:applic}

 We investigate the applications of Theorem \ref{thm: main} to (a) the Regge calculus, and (b) Computer Graphics and related fields.

\subsection{Regge calculus} Looking for the applications of the results above in the context of Regge calculus \cite{reg},  is only natural, given that the original motivation of \cite{cms} (and its precursory, more Physics oriented, \cite{cms1}) stemmed precisely from therein. %
Indeed, as emphasized in \cite{cms}, Gauss (scalar) curvature $K \equiv R^2$ and mean curvature $H \equiv H^2$ are relevant to the Hilbert action principle, thence to the derivation of the Einstein field equations.

\subsubsection{An immediate consequence} \label{subsec:wormholes}

As we have already noted above, our result represents a ``reverse version'' of a result of \cite{cms}.

Its relevance to (Discrete) General Relativity is accentuated by contemporary works such  as \cite{Sm} -- see also the discussion in Section \ref{sec:def-fat} above. Indeed, our main result \ref{thm: main} has as a (rather simple) consequence the possibility of extending the intrinsic the intrinsic differential structure 
from the surface singularity to its interior, that is, starting from a horizon of prescribed (and arbitrary) geometry, constructing an asymptotically flat structure.
This is possible not only for the immediate (rather trivial) case of smooth
singularities, e.g. {\it Wheeler wormholes} \cite{Wh1} (i.e. such that $S \simeq \mathbb{S}^2$, where $S \subset \partial M^3$) that constituted an important motivational phenomenon in Regge's original paper \cite{reg}, but also for the more interesting case of {\it Wheeler
Foam} \cite{Wh}, i.e. non-smooth singularities (for which topology and
curvature change with scale).\footnote{An amusing (almost ``Sci-Fi") consequence of this result is the fact that it allows for a continuous ``gluing'' of the geometry of the given ``Universe'', across the black-hole singularity (i.e. the common boundary), to the ``Alternate Universe'', without any distortion (hence without observable changes for an observer), except on a arbitrarily small, symmetric tubular neighbourhood of the singularity.}

\subsubsection{Purely metric Regge calculus}

The Physical motivation of the study of the curvatures convergence problem raises a few natural, interrelated, questions:

\begin{quest}
The $PL$ spaces are still not the discrete metric spaces (lattices) sought for in
quantum field theory. Can one discard this restriction?
\end{quest}

The answer turns out to be affirmative, as we indicate below, if one discards the angular defect approach in favor of working with metric curvatures

\begin{quest}
Since Regge's drive was to find a purely metric (discrete)
formulation of Gravity, the presence of angles in
the Lipschitz-Killing curvatures is a bit ``unesthetic'', as already stressed above.
Hence: can one
(non-trivially) formulate Theorem \ref{thm: main} (and its consequences) solely in
metric terms?
\end{quest}

Again, the answer to this question, as for the previous one, is positive and it rests upon the following formula:

%
\begin{equation}  \label{eq:sectional}
R^j(M^n) =
\frac{1}{{\rm Area}(\mathbb{S}^{n-j-1})}\int_{M^{n-1}}S_{n-j-1}(k_1(x),k_2(x),\ldots,k_{n-1}(x))d\mathcal{H}^{n-1}\,,
\end{equation}
where $M^{n-1} = \partial (M^n$), $d\mathcal{H}^{n-1}$ denotes the $(n-1)$-dimensional Hausdorff measure, and where the {\it symmetric
functions} $S_{j}$ are defined by:
\begin{equation}  \label{eq:sectional1}
S_{j}\big(k_1(x),k_2(x),\ldots,k_{j-1}(x)\big) = \sum_{1 \leq
k_{i_1} \leq k_{i_k} \leq j-1}k_{i_1}(x) \cdots k_{i_k}(x)\,,
\end{equation}
$k_1(x),k_2(x),\ldots,k_{n-1}(x)$ being the principal curvatures  -- see \cite{Za}. (See also \cite{cms2}, Formula (37).)

The formula above shows why the Lipschitz-Killing curvatures are also called the {\it total mean curvatures} (and the ``$S_j$''-s are called the {\it mean curvatures (of order $j$)}. It also suggests a quite direct method of obtaining a local (point-wise) version. However, we can do better than adopting this somewhat naive approach and, moreover, extending the definition to a wide range of geometric objects, by making appeal to methods and results of Geometric Measure Theory, which we succinctly present here (see \cite{Za} and, for more details, \cite{Za1}):

\begin{defn}
Let $X \subset \mathbb{R}^d$. The {\it reach} of $X$ is defined as:
\begin{equation}
{\rm reach}(X) = \sup\{r > 0\,|\, 
\forall  y \in X_r\,, \exists!\; x \in X {\rm nearest\; to\;} y \}\,,
\end{equation}
where $X_r$ denotes the $r$-neighbourhood of $X$.

The set of all sets (in $\mathbb{R}^d$) of positive reach is denoted by $PR$.
\end{defn}

We can now introduce the notation for the Lipschitz-Killing curvature in the desired generalized context:

Let $X \in PR$, let $\varepsilon < {\rm reach}(X)$ and let $B \subset \mathbb{R}^d$ be a Borel set. We denote by $C_k(X,B)$ the {\it Lipschitz-Killing curvature measures}\footnote{They are signed Radon measures on $B$}. We do not define them formally, but rather present below the characterization theorem we are interested in for our ends. However, we first have to make a few observations as well as introduce a little bit more notations.

First, let us note that, small enough $\varepsilon > 0$, $\partial X_r$ is a 
$\mathcal{C}^{1,1}$-hypersurface.\footnote{For more details see, for instance, \cite{lse} and the bibliography therein.}
          Therefore, they admit principal curvatures (in the classical sense) $k^\varepsilon_i(x + \mathbf{n})$ at almost any point
$p = x + \mathbf{n}$, where $\mathbf{n}$ denotes the normal unit vector (at $x$).
          Define the generalized principal curvatures by: $k_i(\varepsilon,\mathbf{n}) = k^\varepsilon_i(x + \mathbf{n})$.
          Then $k_i(\varepsilon,\mathbf{n})$ exist $\mathcal{H}^{n-1}$-a.a. $(x,\mathbf{n})$. (Here $\mathcal{H}$ denotes the Hausdorff measure.)\footnote{See also \cite{lse} for the application of this approach in Computer Graphics.}

Recall also the following

\begin{defn}
The ({\it unit}) {\it normal bundle} of $X$ is defined as:
\begin{equation}
{\rm nor}(X) = \{(x,\mathbf{n}) \in \partial X \times \mathbb{S}^{d-1}\,|\,\mathbf{n} \in {\rm Nor}(X,x)\},
\end{equation}
where ${\rm Nor}(X,x) = \{\mathbf{n} \in S_{d-1}\,|\, <\mathbf{n},v> \leq 0, v \in {\rm Tan}(X,x))\}$ is the {\it normal cone} (to $X$ at the point $x \in T$), dual to the {\it tangent cone} (to $X$ at the point $x \in T$) -- see, \cite{Fe}.
\end{defn}

\begin{thm}
With the notations above
\begin{equation} \label{eq:KLvsKi}
C_k(X,B) = \int_{{\rm nor}(X)}\mathbf{1}_B\prod_{i=1}^{d-1}\frac{1}{\sqrt{1+\kappa_i(x,\mathbf{n})^2}}S_{d-1-k}(\kappa_1(x,\mathbf{n}),\cdots,\kappa_{d-1}(x,\mathbf{n}))d\mathcal{H}^{d-1}(x,\mathbf{n})\,.
\end{equation}
\end{thm}

Now, our suggested approach is, evidently, to use formula (\ref{eq:sectional}) to express the curvature at a vertex by principal curvatures (which is easy -- and natural -- to determine for a piecewise-flat surface). Use then results of 
\cite{BM}, Section 10, to show that piecewise-flat, metric curvatures (Menger and/or Haantjes) can approximate well (in fact: as well as desired) smooth (classical) curvature\footnote{and also torsion} (and, of course, arc length of curves). In consequence one can compute the Gauss and mean curvatures of piecewise-flat surfaces, with a clear and well established importance in Graphics, etc. (See \cite{lse} for a discussion and  \cite{SA} for an implementation and some numerical results.)

However, these are only approximations with limited 
convergence properties -- see \cite{SA}, \cite{Sa04} for the metric curvatures aspect and also \cite{Su++}, \cite{D++}, amongst others, for the more general problem of approximating principal curvatures. The problem is that, in order to ensure convergence for the Cheeger et al. process, the triangulation has not only to converge {\it in mesh} to $0$, it also has to remain fat. However, to ensure a ``good sampling'' of the directions on a surface (so to ensure good approximation of the principal curvatures), one necessarily has to produce samplings whose edges directions are arbitrarily dense in the tangent plane (cone) at a given point (vertex), in contradiction with the previous fatness constraint.
     (Note, however, that this problem does not exist if one is willing 
to be content with approximate results, correct up to a predetermined error. One should also keep in mind Remarks \ref{rem:big-fat} and \ref{rem:Fu-fat} regarding Fu's work, as relevant 
to this aspect.)

     Moreover, given the fact (discussed above), that one can not increase as desired the number of triangles adjacent to a vertex, there exist only a very limited number of directions to ``choose'' from. 
     Thus one is quite restricted in adding new directions (thus to better ``sample'' the surface, so to say) without negatively affecting the fatness of the triangulation.\footnote{In practice (Graphics, etc.) even angles $\pi/12$ are already problematically small!...}

     How to ``get rid'' of this problem?
     \begin{itemize}

     \item ``Mix'' the angles in the manner described in detail in \cite{s06}, to obtain angles whose measure is close to their mean (i.e. $\pi/{\rm deg}(v)$ -- where ${\rm deg}(v)$ denotes the number of triangles adjacent to the vertex $v$. (This ``trick'' works if a  smoothness condition (albeit, minimal) is imposed on the manifold -- it certainly holds for triangulated surfaces.)

     \item Add directions by considering the $PL$-{\it quasi-geodesic} 
and ``normalize'' (by projection on the normal plane) -- see \cite{LocGlob1}, \cite{SA} and the references therein for details and some numerical results.

     \item The approximation of principal curvatures being, as already noted above,
     notoriously difficult -- see \cite{Ma} \cite{Su++}, \cite{Po++}, \cite{D++} -- one seeks other, perhaps less direct strategies.
     Such a method does, in indeed exist, embodied by the generalized principal curvatures, as we have detailed above. Passing to smooth surfaces allows for the use of a wide scale of well developed and finely honed methods of Graphics and related fields. Furthermore, it clearly compensates for its departure from the given (data) set of discrete/geometric $PL$ object by its generality whence its applicative potential in a very general setting.


     \end{itemize}

As we have seen above, we can compute the principal curvatures of $S$ via those of the smoother surface $S_\varepsilon$, at least up to some infinitesimal distortion. However, to determine the full curvature tensor, in the case of higher dimensional manifolds, suffices to determine the Gauss curvature of all the 2-sections (see observation above). So, for the full reconstruction of the curvature tensor it is not necessary to determine the principal curvatures, suffices to find the sectional curvatures. We shall show shortly how we this can be done.

\begin{rem}
It is natural to ask the question whether it is possible to compute the Lipschitz-Killing curvatures, starting from the defining Formula (\ref{eq:Rj-smooth}), that is if one can compute the necessary curvature 2-forms and connection 1-forms. The answer seems to be positive, even though, till recently this was only a mainly theoretical possibility (see \cite{So}). However, after the appearance of the computational exterior differential calculus, introduced by Gu \cite{Gu} and Gu and Yau \cite{GY}, and embraced and developed since then by many others, this approach appears quite feasible, at least in dimensions 2 and 3.

Another related question is whether it is possible to determine other important curvatures (e.g. Ricci and scalar) as well as the Lipschitz-Killing curvatures, via the sectional curvatures, without appeal to principal curvatures (and Formula (\ref{eq:sectional}) or (\ref{eq:KLvsKi})\footnote{or, (\ref{eq:Rj}), for that matter.}).
We discuss this problem separately in a related paper \cite{Sa11d} 
and present a solution for the case of Ricci curvature of $PL$ manifolds.

\end{rem}

\subsection{Computer Graphics} \label{subsec:graphics}
Polygonal/polyhedral (mainly triangular/tetrahedral) meshes are the basic representations of geometry, employed in a plethora of related  fields, such as Computer Vision, Image Processing, Computer Graphics, Geometric Modeling and Manufacturing. Curvature analysis of this type of data sets plays a major role in a variety of applications, such as reconstruction, segmentation and recognition and non photorealistic
rendering (see the bibliography included in \cite{lse} for some of the vast -- and ever developing -- literature on the subject).

Since, as already mentioned above, $R^0 \equiv$ {\rm volume}, $R^2
\equiv$ {\rm scalar curvature}
and $H^1 \equiv$ {area boundary},
$H^2 \equiv$ mean curvature for inward normal,
etc., there exists a relationship between the main subject of this paper and the fields mentioned above.
Indeed, one such connection was already mentioned in Remark \ref{rem:extension} above.
However the connection is deeper and less trivial than this. Indeed, in Computer Graphics, Computer Aided Geometric Design, etc. it has become customary 
lately to compute so called {\it volumetric curvatures} (see, e.g. \cite{SER}), and Graphics (see \cite{GY1} and the bibliography therein). This amounts, in fact, to the computation of the curvatures (Gauss and mean) of surfaces evolving in time. We take this opportunity to note that, while surely this approach has merit and uses, a proper ``volumetric'', i.e. 3-dimensional curvature (measures), would entail the computation of a namely {\it sectional}, {\it scalar} and {\it Ricci} curvatures. (It seems, indeed, that the last one deserves special attention, at least in Image Processing -- see \cite{SAWZ1}.)
In all fairness, we should add that, for Computer Graphics, where usually the data is already embedded in $\mathbb{R}^3$, thus endowed with the Euclidean (flat) geometry of the ambient space, such computations are, therefore, 
rather meaningless.

Clearly, the metric approach to the computation of $PL$ (and polyhedral) manifolds considered above is more relevant in practice for the applicative fields considered above. This is particularly true for 2-dimensional manifolds, with further emphasis on the types of surface (square grids) traditionally employed in Imaging, since, as we have already noted, the notion of Ricci curvature for the dual complex and the one for the original (given) complex, have the same geometric significance and even coincide, perhaps up to a constant.

We conclude this section by noting that estimation of curvatures, mainly of mean curvature, is also important (mostly via the {\it Cahn-Hilliard equation}) in physically motivated applications -- see, e.g. \cite{Ra}, amongst many others. 


\section{Quasiconformal Mappings}


As a consequence of Theorem \ref{thm:eu}, we have proved in \cite{s2} the following existence result

\begin{thm} [\cite{s2}]
Let $M^n$, Riemannian manifold, where $M$ satisfies the conditions in the statement of Theorem \ref{thm:eu}.
Then there exists a non-constant quasimeromorphic mapping $f:M^n
\rightarrow \widehat{\mathbb{R}^n}$.
\end{thm}
(Here $\widehat{\mathbb{R}^n} = \mathbb{R}^n \cup \{\infty\}$ is identified with $\mathbb{S}^n$ with spherical metric.)

Since the proof essentially requires the construction of a ``chess-board'' fat triangulation (followed by the alternate quasiconformal mapping of the ``black'' and ``white''' simplices to the interior, respective exterior of the standard simplex in $\mathbb{R}^n$), it follows that Theorem \ref{thm:cms} and our own Therems \ref{thm:eu} and \ref{thm: main} show that spaces that admit
``good'' curvature convergence in secant approximation are
{\it geometric branched covers} of $\mathbb{S}^n$.

Recall that

\begin{defn} \label{def:QR}
Let $(M,d), (N,\rho)$ be metric spaces and let $f:(M,d) \rightarrow (N^n,\rho)$ be a continuous function. Then $f$ is called
\begin{enumerate}
\item {\it quasiregular} (or, more precisely, $K$-quasiregular) iff there exists $1 \leq K \leq \infty$, such that, for any $x \in M$, the following holds:
\begin{equation}
H(f,x) = \limsup_{r \rightarrow 0}\frac{\sup\{\rho(f(x),f(y))\;|\;d(x,y) = r\}}{\inf\{\rho(f(x),f(y))\;|\;d(x,y) = r\}} \leq K;
\end{equation}
\item {\it quasiconformal} iff it is a quasiregular homeomorphism;
\item {\it quasimeromorphic} iff $(N,\rho)$ is the unit sphere $S^n$ (equipped with standard metric).
\end{enumerate}
$H(f,x)$ is called the {\it linear dilatation} of $f$ (at $x$).
\end{defn}

Obviously, the linear dilatation is a measure of the eccentricity of the image of infinitesimal balls. Therefore, (at least if one restricts oneself to Riemannian manifolds) quasiconformal mappings can be characterized as being precisely those maps that
\begin{itemize}
\item map 
infinitesimal balls into infinitesimal ellipsoids (of bounded eccentricity);

\item map 
almost balls into almost ellipsoids;

\item  distort infinitesimal spheres by a constant factor.
\end{itemize}
In fact, if one considers (in the Riemannian manifold setting) the linear mapping $f':\mathbb{R}^n \rightarrow \mathbb{R}^n$, then $f'(B^n) = E(f')$ is an ellipsoid of semi-axes $a_1 \geq a_2 \geq \cdots \geq a_n$\footnote{and equal to the square roots of the eigenvalue of the adjoint mapping of $f'$} and the characterizations above follow. Not only this, but, in fact,
\begin{equation} \label{eq:H(a1,an)}
H(A) = \frac{a_1}{a_n}.
\end{equation}
Moreover, quasiconformal mappings
\begin{itemize}
\item distort local distances by a fixed amount;
\item preserve approximative shape.\footnote{but in these cases the characterization is not sharp, the proper class of functions characterized by this property being the so called {\it quasisymmetric mapings} -- see, e.g. \cite{Gr-carte}, p. 418 ff.}
\end{itemize}

\begin{rem}
There exist two other definitions of quasiconformality (for mappings between Riemannian manifolds of the same dimensionality), but we have chosen the one above -- the so called {\it metric definition} -- due to its simplicity, naturalness in our context (see below) and the fact that it makes sense for any metric space. On the other hand, if one wishes to prove even the simplest, intuitive geometric properties (like the one mentioned below),
a delicate interplay of all of the three definitions is needed.
\end{rem}

Of course, one naturally asks whether the ``quasiconformal'' in Definition \ref{def:QR} above implies, indeed, as the name suggests, that quasiconformal mappings ``almost'' preserve angles (given that {\it conformal} mappings do). The answer is, as expected, positive -- see \cite{Ag}.

Given this last intuitive characterization of quasiconformality, and Proposition \ref{eq:fat-CMS}, in conjunction with Theorem \ref{thm:eu} and Theorem \ref{thm: main} readily implies

\begin{prop} \label{prop:CQ}
Let $M_1^n,M_2^n$ be two (connected) manifolds, topologically as in the statement of Theorem \ref{thm:eu}, and let $\varphi: M_1^n \rightarrow M^n_2$ be a $K$-quasiconformal mapping.
Then
\begin{enumerate}
\item If $M_1^n,M_2^n$ are $PL$ manifolds, then
\begin{equation} \label{eq:qc-tensor}
R^j_2 = C(K,n,j)R^j_1\,,
\end{equation}
where $R^j_1, R^j_2$ denote the Lipschitz-Killing curvatures of $M_1^n, M_2^n$, respectively.
\item If $M_1^n,M_2^n$ are smooth manifolds manifolds, then
\begin{equation} \label{eq:qc-tensor}
R^j_2 = (1 + \varepsilon_0)C(K,n,j)R^j_1\,,
\end{equation}
in the sense of measures, for some arbitrarily small $\varepsilon_1 > 0$.
\end{enumerate}
\end{prop}

Alternatively, in view of the remarks above, it is easily to see that Formula (\ref{eq:sectional}), in conjunction with the definition of quasiconformality, also readily implies the result above. In this approach, one regards the relevant edges of a fine enough triangulation both as the principal vectors and as semi-axes of an ``infinitesimal' ellipsoid.

Let us note ``en passant'' that Theorem \ref{thm:aRs} and the technique in \cite{s2}, immediately imply the following

\begin{cor}
Let $(M,d)$ be a connected almost Riemannian manifold, topologically as in the statement of Theorem \ref{thm:eu}. Then there exists a non-constant quasimeromorphic mapping $f:M
\rightarrow \widehat{\mathbb{R}^n}$.
\end{cor}

Proposition \ref{prop:CQ} above being, clearly, a direct consequence of Theorem \ref{thm:eu}, it could have been placed, based on this criterion, in the first section. However, given that it concerns primarily mappings, rather than curvature, we have decided to included in a different section of its own, in order to not alter the natural flow of ideas and concepts. Also, not to depart to much from the main framework of the paper, we do not consider here the problem of extending Proposition \ref{prop:CQ} to the case of general quasiregular mappings, as well as discussing other related aspects relating quasiregular mappings and curvature, but we rather postpone them for further study 
-- see \cite{Sa11a}.

We conclude, therefore, this paper with the following Remark and the ensuing Question:

\begin{rem}
Since all the theorems for geometric branched coverings of
$\mathbb{S}^n$ were obtained via the Alexander Trick, i.e. by
constructing fat chessboard triangulation one is easily conducted to
the following:
\end{rem}

\begin{quest}
Does $M^n$ admit a $qr$-mapping on $\mathbb{S}^n$ iff it admits
``good'' curvature convergence in secant approximation?
\end{quest}




\end{document}